\documentclass{conm-p-l}

\theoremstyle{definition}

\theoremstyle{remark}

\numberwithin{equation}{section}

\usepackage{graphicx,amsthm}
\usepackage{verbatim,amsmath,amscd,amssymb,color}
\input xy
\xyoption{all}
\begin{document}
\renewcommand{\labelenumi}{$($\roman{enumi}$)$}
\renewcommand{\labelenumii}{$(${\rm \alph{enumii}}$)$}
\font\germ=eufm10
\newcommand{\cI}{{\mathcal I}}
\newcommand{\cA}{{\mathcal A}}
\newcommand{\cB}{{\mathcal B}}
\newcommand{\cC}{{\mathcal C}}
\newcommand{\cD}{{\mathcal D}}
\newcommand{\cE}{{\mathcal E}}
\newcommand{\cF}{{\mathcal F}}
\newcommand{\cG}{{\mathcal G}}
\newcommand{\cH}{{\mathcal H}}
\newcommand{\cK}{{\mathcal K}}
\newcommand{\cL}{{\mathcal L}}
\newcommand{\cM}{{\mathcal M}}
\newcommand{\cN}{{\mathcal N}}
\newcommand{\cO}{{\mathcal O}}
\newcommand{\cR}{{\mathcal R}}
\newcommand{\cS}{{\mathcal S}}
\newcommand{\cV}{{\mathcal V}}
\newcommand{\fra}{\mathfrak a}
\newcommand{\frb}{\mathfrak b}
\newcommand{\frc}{\mathfrak c}
\newcommand{\frd}{\mathfrak d}
\newcommand{\fre}{\mathfrak e}
\newcommand{\frf}{\mathfrak f}
\newcommand{\frg}{\mathfrak g}
\newcommand{\frh}{\mathfrak h}
\newcommand{\fri}{\mathfrak i}
\newcommand{\frj}{\mathfrak j}
\newcommand{\frk}{\mathfrak k}
\newcommand{\frI}{\mathfrak I}
\newcommand{\fm}{\mathfrak m}
\newcommand{\frn}{\mathfrak n}
\newcommand{\frp}{\mathfrak p}
\newcommand{\fq}{\mathfrak q}
\newcommand{\frr}{\mathfrak r}
\newcommand{\frs}{\mathfrak s}
\newcommand{\frt}{\mathfrak t}
\newcommand{\fru}{\mathfrak u}
\newcommand{\frA}{\mathfrak A}
\newcommand{\frB}{\mathfrak B}
\newcommand{\frF}{\mathfrak F}
\newcommand{\frG}{\mathfrak G}
\newcommand{\frH}{\mathfrak H}
\newcommand{\frJ}{\mathfrak J}
\newcommand{\frN}{\mathfrak N}
\newcommand{\frP}{\mathfrak P}
\newcommand{\frT}{\mathfrak T}
\newcommand{\frU}{\mathfrak U}
\newcommand{\frV}{\mathfrak V}
\newcommand{\frX}{\mathfrak X}
\newcommand{\frY}{\mathfrak Y}
\newcommand{\frZ}{\mathfrak Z}
\newcommand{\rA}{\mathrm{A}}
\newcommand{\rC}{\mathrm{C}}
\newcommand{\rd}{\mathrm{d}}
\newcommand{\rB}{\mathrm{B}}
\newcommand{\rD}{\mathrm{D}}
\newcommand{\rE}{\mathrm{E}}
\newcommand{\rH}{\mathrm{H}}
\newcommand{\rK}{\mathrm{K}}
\newcommand{\rL}{\mathrm{L}}
\newcommand{\rM}{\mathrm{M}}
\newcommand{\rN}{\mathrm{N}}
\newcommand{\rR}{\mathrm{R}}
\newcommand{\rT}{\mathrm{T}}
\newcommand{\rZ}{\mathrm{Z}}
\newcommand{\bbA}{\mathbb A}
\newcommand{\bbB}{\mathbb B}
\newcommand{\bbC}{\mathbb C}
\newcommand{\bbG}{\mathbb G}
\newcommand{\bbF}{\mathbb F}
\newcommand{\bbH}{\mathbb H}
\newcommand{\bbP}{\mathbb P}
\newcommand{\bbN}{\mathbb N}
\newcommand{\bbQ}{\mathbb Q}
\newcommand{\bbR}{\mathbb R}
\newcommand{\bbV}{\mathbb V}
\newcommand{\bbZ}{\mathbb Z}
\newcommand{\adj}{\operatorname{adj}}
\newcommand{\Ad}{\mathrm{Ad}}
\newcommand{\Ann}{\mathrm{Ann}}
\newcommand{\rcris}{\mathrm{cris}}
\newcommand{\ch}{\mathrm{ch}}
\newcommand{\coker}{\mathrm{coker}}
\newcommand{\diag}{\mathrm{diag}}
\newcommand{\Diff}{\mathrm{Diff}}
\newcommand{\Dist}{\mathrm{Dist}}
\newcommand{\rDR}{\mathrm{DR}}
\newcommand{\ev}{\mathrm{ev}}
\newcommand{\Ext}{\mathrm{Ext}}
\newcommand{\cExt}{\mathcal{E}xt}
\newcommand{\fin}{\mathrm{fin}}
\newcommand{\Frac}{\mathrm{Frac}}
\newcommand{\GL}{\mathrm{GL}}
\newcommand{\Hom}{\mathrm{Hom}}
\newcommand{\hd}{\mathrm{hd}}
\newcommand{\rht}{\mathrm{ht}}
\newcommand{\id}{\mathrm{id}}
\newcommand{\im}{\mathrm{im}}
\newcommand{\inc}{\mathrm{inc}}
\newcommand{\ind}{\mathrm{ind}}
\newcommand{\coind}{\mathrm{coind}}
\newcommand{\Lie}{\mathrm{Lie}}
\newcommand{\Max}{\mathrm{Max}}
\newcommand{\mult}{\mathrm{mult}}
\newcommand{\op}{\mathrm{op}}
\newcommand{\ord}{\mathrm{ord}}
\newcommand{\pt}{\mathrm{pt}}
\newcommand{\qt}{\mathrm{qt}}
\newcommand{\rad}{\mathrm{rad}}
\newcommand{\res}{\mathrm{res}}
\newcommand{\rgt}{\mathrm{rgt}}
\newcommand{\rk}{\mathrm{rk}}
\newcommand{\SL}{\mathrm{SL}}
\newcommand{\soc}{\mathrm{soc}}
\newcommand{\Spec}{\mathrm{Spec}}
\newcommand{\St}{\mathrm{St}}
\newcommand{\supp}{\mathrm{supp}}
\newcommand{\Tor}{\mathrm{Tor}}
\newcommand{\Tr}{\mathrm{Tr}}
\newcommand{\wt}{\mathrm{wt}}
\newcommand{\Ab}{\mathbf{Ab}}
\newcommand{\Alg}{\mathbf{Alg}}
\newcommand{\Grp}{\mathbf{Grp}}
\newcommand{\Mod}{\mathbf{Mod}}
\newcommand{\Sch}{\mathbf{Sch}}\newcommand{\bfmod}{{\bf mod}}
\newcommand{\Qc}{\mathbf{Qc}}
\newcommand{\Rng}{\mathbf{Rng}}
\newcommand{\Top}{\mathbf{Top}}
\newcommand{\Var}{\mathbf{Var}}
\newcommand{\gromega}{\langle\omega\rangle}
\newcommand{\lbr}{\begin{bmatrix}}
\newcommand{\rbr}{\end{bmatrix}}
\newcommand{\forb}{\bigcirc\kern-2.8ex \because}
\newcommand{\forbb}{\bigcirc\kern-3.0ex \because}
\newcommand{\forbbb}{\bigcirc\kern-3.1ex \because}
\newcommand{\cd}{commutative diagram }
\newcommand{\SpS}{spectral sequence}
\newcommand\C{\mathbb C}
\newcommand\hh{{\hat{H}}}
\newcommand\eh{{\hat{E}}}
\newcommand\F{\mathbb F}
\newcommand\fh{{\hat{F}}}
\def\ge{\frg}
\def\AA{{\mathcal A}}
\def\al{\alpha}
\def\bq{B_q(\ge)}
\def\bqm{B_q^-(\ge)}
\def\bqz{B_q^0(\ge)}
\def\bqp{B_q^+(\ge)}
\def\beneme{\begin{enumerate}}
\def\beq{\begin{equation}}
\def\beqn{\begin{eqnarray}}
\def\beqnn{\begin{eqnarray*}}
\def\bigsl{{\hbox{\fontD \char'54}}}
\def\bbra#1,#2,#3{\left\{\begin{array}{c}\hspace{-5pt}
#1;#2\\ \hspace{-5pt}#3\end{array}\hspace{-5pt}\right\}}
\def\cd{\cdots}
\def\CC{\mathbb{C}}
\def\CBL{\cB_L(\TY(B,1,n+1))}
\def\CBM{\cB_M(\TY(B,1,n+1))}
\def\CVL{\cV_L(\TY(D,1,n+1))}
\def\CVM{\cV_M(\TY(D,1,n+1))}
\def\ddd{\hbox{\germ D}}
\def\del{\delta}
\def\Del{\Delta}
\def\Delr{\Delta^{(r)}}
\def\Dell{\Delta^{(l)}}
\def\Delb{\Delta^{(b)}}
\def\Deli{\Delta^{(i)}}
\def\Delre{\Delta^{\rm re}}
\def\ei{e_i}
\def\eit{\tilde{e}_i}
\def\eneme{\end{enumerate}}
\def\ep{\epsilon}
\def\eeq{\end{equation}}
\def\eeqn{\end{eqnarray}}
\def\eeqnn{\end{eqnarray*}}
\def\fit{\tilde{f}_i}
\def\FF{{\rm F}}
\def\ft{\tilde{f}}
\def\gau#1,#2{\left[\begin{array}{c}\hspace{-5pt}#1\\
\hspace{-5pt}#2\end{array}\hspace{-5pt}\right]}
\def\gl{\hbox{\germ gl}}
\def\hom{{\hbox{Hom}}}
\def\ify{\infty}
\def\io{\iota}
\def\kp{k^{(+)}}
\def\km{k^{(-)}}
\def\llra{\relbar\joinrel\relbar\joinrel\relbar\joinrel\rightarrow}
\def\lan{\langle}
\def\lar{\longrightarrow}
\def\max{{\rm max}}
\def\lm{\lambda}
\def\Lm{\Lambda}
\def\mapright#1{\smash{\mathop{\longrightarrow}\limits^{#1}}}
\def\Mapright#1{\smash{\mathop{\Longrightarrow}\limits^{#1}}}
\def\mm{{\bf{\rm m}}}
\def\nd{\noindent}
\def\nn{\nonumber}
\def\nnn{\hbox{\germ n}}
\def\catob{{\mathcal O}(B)}
\def\oint{{\mathcal O}_{\rm int}(\ge)}
\def\ot{\otimes}
\def\op{\oplus}
\def\opi{\ovl\pi_{\lm}}
\def\osigma{\ovl\sigma}
\def\ovl{\overline}
\def\plm{\Psi^{(\lm)}_{\io}}
\def\qq{\qquad}
\def\q{\quad}
\def\qed{\hfill\framebox[2mm]{}}
\def\QQ{\mathbb Q}
\def\qi{q_i}
\def\qii{q_i^{-1}}
\def\ra{\rightarrow}
\def\ran{\rangle}
\def\rlm{r_{\lm}}
\def\ssl{\hbox{\germ sl}}
\def\slh{\widehat{\ssl_2}}
\def\ti{t_i}
\def\tii{t_i^{-1}}
\def\til{\tilde}
\def\tm{\times}
\def\tt{\frt}
\def\TY(#1,#2,#3){#1^{(#2)}_{#3}}
\def\ua{U_{\AA}}
\def\ue{U_{\vep}}
\def\uq{U_q(\ge)}
\def\uqp{U'_q(\ge)}
\def\ufin{U^{\rm fin}_{\vep}}
\def\ufinp{(U^{\rm fin}_{\vep})^+}
\def\ufinm{(U^{\rm fin}_{\vep})^-}
\def\ufinz{(U^{\rm fin}_{\vep})^0}
\def\uqm{U^-_q(\ge)}
\def\uqmq{{U^-_q(\ge)}_{\bf Q}}
\def\uqpm{U^{\pm}_q(\ge)}
\def\uqq{U_{\bf Q}^-(\ge)}
\def\uqz{U^-_{\bf Z}(\ge)}
\def\ures{U^{\rm res}_{\AA}}
\def\urese{U^{\rm res}_{\vep}}
\def\uresez{U^{\rm res}_{\vep,\ZZ}}
\def\util{\widetilde\uq}
\def\uup{U^{\geq}}
\def\ulow{U^{\leq}}
\def\bup{B^{\geq}}
\def\blow{\ovl B^{\leq}}
\def\vep{\varepsilon}
\def\vp{\varphi}
\def\vpi{\varphi^{-1}}
\def\VV{{\mathcal V}}
\def\xii{\xi^{(i)}}
\def\Xiioi{\Xi_{\io}^{(i)}}
\def\W1{W(\varpi_1)}
\def\WW{{\mathcal W}}
\def\wt{{\rm wt}}
\def\wtil{\widetilde}
\def\what{\widehat}
\def\wpi{\widehat\pi_{\lm}}
\def\ZZ{\mathbb Z}

\def\m@th{\mathsurround=0pt}
\def\fsquare(#1,#2){
\hbox{\vrule$\hskip-0.4pt\vcenter to #1{\normalbaselines\m@th
\hrule\vfil\hbox to #1{\hfill$\scriptstyle #2$\hfill}\vfil\hrule}$\hskip-0.4pt
\vrule}}

\newtheorem{thm}{Theorem}[section]
\newtheorem{pro}[thm]{Proposition}
\newtheorem{lem}[thm]{Lemma}
\newtheorem{ex}[thm]{Example}
\newtheorem{cor}[thm]{Corollary}
\newtheorem{conj}[thm]{Conjecture}
\theoremstyle{definition}
\newtheorem{df}[thm]{Definition}

\newcommand{\cmt}{\marginpar}
\newcommand{\seteq}{\mathbin{:=}}
\newcommand{\cl}{\colon}
\newcommand{\be}{\begin{enumerate}}
\newcommand{\ee}{\end{enumerate}}
\newcommand{\bnum}{\be[{\rm (i)}]}
\newcommand{\enum}{\ee}
\newcommand{\ro}{{\rm(}}
\newcommand{\rf}{{\rm)}}
\newcommand{\set}[2]{\left\{#1\,\vert\,#2\right\}}
\newcommand{\sbigoplus}{{\mbox{\small{$\bigoplus$}}}}
\newcommand{\ba}{\begin{array}}
\newcommand{\ea}{\end{array}}
\newcommand{\on}{\operatorname}
\newcommand{\eq}{\begin{eqnarray}}
\newcommand{\eneq}{\end{eqnarray}}
\newcommand{\hs}{\hspace*}

\title{Affine Geometric Crystal of type $G^{(1)}_2$}

\author{Toshiki N\textsc{akashima}}
\address{Department of Mathematics, 
Sophia University, Kioicho 7-1, Chiyoda-ku, Tokyo 102-8554,
Japan}
\email{toshiki@mm.sophia.ac.jp}
\thanks{supported in part by JSPS Grants 
in Aid for Scientific Research.}

\subjclass{Primary 17B37; 17B67; Secondary 22E65; 14M15}
\date{}

\dedicatory{In honor of professor James Lepowsky and 
professor Robert L.Wilson.}

\keywords{geometric crystal, tropicalization, 
ultra-discretization}

\begin{abstract}
We shall realize certain affine geometric 
crystal of type $G^{(1)}_2$ explicitly in the fundamental 
representation $W(\varpi_1)$. 
Its explicit form is rather complicated but still keeps
a positive structure.
\end{abstract}

\maketitle
\renewcommand{\thesection}{\arabic{section}}
\section{Introduction}
\setcounter{equation}{0}
\renewcommand{\theequation}{\thesection.\arabic{equation}}

Geometric crystal is an object defined over certain algebraic
(or ind-)variety which holds an analogous structure to Kashiwara's crystal
(\cite{BK},\cite{N}). Precisely, for a fixed Cartan data
$(A,\{\al_i\}_{i\in I},\{h_i\}_{\i\in I})$, 
a geometric crystal consists of
an ind-variety $X$ over the complex number $\bbC$, a rational
$\bbC^\times$-action $e_i:\bbC^\times\times X\longrightarrow X$ and 
rational functions $\gamma_i,\vep_i:X\longrightarrow \bbC$ $(i\in I)$,
which satisfy certain conditions (see Definition \ref{def-gc}).
It has many similarity to the theory of crystals, {\it e.g.,}
some product structure, Weyl group actions, R-matrices, 
{\it etc}. Furthermore, there is a more direct correspondence
between geometric crystals and crystals, called tropicalization/
ultra-discretization procedure (see \S2).

In \cite{KNO}, we presented certain conjecture. 
In order to mention it precisely, we need to prepare the following:
Let $G$ (resp.~$\ge$) be the affine 
Kac-Moody group 
(resp.~algebra) associated with a 
generalized Cartan matrix $A=(a_{ij})_{i,j\in I}$. Let
 $B^\pm$ be the Borel subgroup and $T$ the 
maximal torus. Set $y_i(c)\seteq\exp(cf_i)$,
and let $\al_i^\vee(c)\in T$ be the image of $c\in\bbC^\times$
by the group morphism $\bbC^\times\to T$ induced by
the simple coroot $\alpha_i^\vee$ as in \ref{KM}.
We set $Y_i(c)\seteq y_i(c^{-1})\,\al_i^\vee(c)=\al_i^\vee(c)\,y_i(c)$.
Let $W$ (resp.~$\wtil W$) be the Weyl group 
(resp.~the extended Weyl group) associated with $\ge$.
The Schubert cell $X_w\seteq BwB/B$ $(w=s_{i_1}\cd s_{i_k}\in W)$ 
is birationally isomorphic to the variety
\[
 B^-_w\seteq\set{Y_{i_1}(x_1)\cd Y_{i_k}(x_k)}%
{x_1,\cd,x_k\in \bbC^\times}\subset B^-,
\]
and $X_w$ has a natural geometric crystal structure
(\cite{BK}, \cite{N}), 

We choose $0\in I$ as in \cite{K0},
and let $\{\varpi_i\}_{i\in I\setminus\{0\}}$
be the set of level $0$ fundamental weights.
Let $W(\varpi_i)$ be the fundamental representation
of $U_q(\ge)$ with $\varpi_i$ as an extremal weight (\cite{K0}).
Let us denote its reduction at $q=1$
by the same notation $W(\varpi_i)$. 
It is a finite-dimensional $\ge$-module.
Note that though the representation $W(\varpi_i)$ is 
irreducible over $\uq$, 
the module $W(\varpi_i)$ 
at $q=1$ for $i\ne1$ is not necessarily an irreducible $\ge$-module.
We set $\bbP(\varpi_i)\seteq
(W(\varpi_i)\setminus\{0\})/\bbC^\times$.

For any $i\in I$, define 
\begin{eqnarray}
&& c_i^\vee\seteq
\mathrm{max }(1,\frac{2}{(\al_i,\al_i)}).
\label{eq:ci}
\end{eqnarray}
Then the translation $t(c^\vee_i\varpi_i)$ 
belongs to $\widetilde W$ (see \cite{KNO}).
For a subset $J$ of $I$, let us denote by
$\ge_J$ the subalgebra of $\ge$ generated by
$\{e_i,f_i\}_{i\in J}$.
For an integral weight $\mu$, define 
$I(\mu)\seteq\set{j\in I}{\lan h_j,\mu\ran\geq0}$.

\begin{conj}[\cite{KNO}]
For any $i\in I$, 
there exist a unique variety $X$ endowed with
a positive $\ge$-geometric crystal structure and a
rational mapping $\pi\cl X\longrightarrow
\bbP(\varpi_i)$ satisfying the following property:
\begin{enumerate}
\item
for an arbitrary extremal vector $u\in W(\varpi_i)_\mu$,
writing the translation
$t(c_i^\vee\mu)$ as $\io w\in 
\wtil W$ with a Dynkin diagram automorphism $\io$
and $w=s_{i_1}\cd s_{i_k}$,
there exists a birational mapping
$\xi\cl B^-_w\longrightarrow X$
such that $\xi$ is a morphism of $\ge_{I(\mu)}$-geometric crystals
and that
the composition
$\pi\circ\xi\cl B^-_w\to \bbP(\varpi_i)$
coincides with
$Y_{i_1}(x_1)\cdots Y_{i_k}(x_k)\mapsto 
Y_{i_1}(x_1)\cd Y_{i_k}(x_k)\ovl u$,
where $\ovl u$ is the line including $u$,
\item
the ultra-discretization of $X$
is isomorphic to the crystal $B_\infty(\varpi_i)$
of the Langlands dual $\ge^L$.
\end{enumerate}
\end{conj}
In \cite{KNO}, we constructed
a positive geometric crystal $\cV(\ge)$ associated 
with the fundamental representation $W(\varpi_1)$ for affine 
Lie algebras $\ge=\TY(A,1,n), \TY(B,1,n),\TY(C,1,n),\TY(D,1,n),$
$\TY(A,2,2n-1),\TY(A,2,2n),\TY(D,2,n+1)$
with this conjecture as a guide.
In that article, we also show that the ultra-discretization 
limit of $\cV(\ge)$ 
is isomorphic to the limit of certain coherent family 
of perfect crystals
for $\ge^L$ the Langlands dual of $\ge$.

In this article, we shall construct such geometric crystal for 
$\ge=\TY(G,1,2)$.
Its explicit form is given in \S5, which is rather complicated 
but we shall see that it is positive, 
which implies that the former half of 
our conjecture is affirmative 
for $\TY(G,1,2)$ and the $i=1$-case. 
Then we obtain its ultra-discretization limit and 
we expect that it is isomorphic to the limit of certain coherent 
family of perfect crystals of type $\TY(D,3,4)$ \cite{KMOY}.

\renewcommand{\thesection}{\arabic{section}}
\section{Geometric crystals}
\setcounter{equation}{0}
\renewcommand{\theequation}{\thesection.\arabic{equation}}

In this section, 
we review Kac-Moody groups and geometric crystals
following 
\cite{PK}, \cite{Ku2}, \cite{BK}
\subsection{Kac-Moody algebras and Kac-Moody groups}
\label{KM}
Fix a symmetrizable generalized Cartan matrix
 $A=(a_{ij})_{i,j\in I}$ with a finite index set $I$.
Let $(\tt,\{\al_i\}_{i\in I},\{\al^\vee_i\}_{i\in I})$ 
be the associated
root data, where ${\tt}$ is a vector space 
over $\bbC$ and
$\{\al_i\}_{i\in I}\subset\tt^*$ and 
$\{\al^\vee_i\}_{i\in I}\subset\tt$
are linearly independent 
satisfying $\al_j(\al^\vee_i)=a_{ij}$.

The Kac-Moody Lie algebra $\ge=\ge(A)$ associated with $A$
is the Lie algebra over $\bbC$ generated by $\tt$, the 
Chevalley generators $e_i$ and $f_i$ $(i\in I)$
with the usual defining relations (\cite{KP},\cite{PK}).
There is the root space decomposition 
$\ge=\bigoplus_{\al\in \tt^*}\ge_{\al}$.
Denote the set of roots by 
$\Delta:=\{\al\in \tt^*|\al\ne0,\,\,\ge_{\al}\ne(0)\}$.
Set $Q=\sum_i\bbZ \al_i$, $Q_+=\sum_i\bbZ_{\geq0} \al_i$,
$Q^\vee:=\sum_i\bbZ \al^\vee_i$
and $\Delta_+:=\Delta\cap Q_+$.
An element of $\Delta_+$ is called 
a {\it positive root}.
Let $P\subset \tt^*$ be a weight lattice such that 
$\bbC\ot P=\tt^*$, whose element is called a
weight.

Define simple reflections $s_i\in{\rm Aut}(\tt)$ $(i\in I)$ by
$s_i(h):=h-\al_i(h)\al^\vee_i$, which generate the Weyl group $W$.
It induces the action of $W$ on $\tt^*$ by
$s_i(\lm):=\lm-\lm(\al^\vee_i)\al_i$.
Set $\Delre:=\{w(\al_i)|w\in W,\,\,i\in I\}$, whose element 
is called a real root.

Let $\ge'$ be the derived Lie algebra 
of $\ge$ and let 
$G$ be the Kac-Moody group associated 
with $\ge'$(\cite{PK}).
Let $U_{\al}:=\exp\ge_{\al}$ $(\al\in \Delre)$
be the one-parameter subgroup of $G$.
The group $G$ is generated by $U_{\al}$ $(\al\in \Delre)$.
Let $U^{\pm}$ be the subgroup generated by $U_{\pm\al}$
($\al\in \Delre_+=\Delre\cap Q_+$), {\it i.e.,}
$U^{\pm}:=\lan U_{\pm\al}|\al\in\Del^{\rm re}_+\ran$.

For any $i\in I$, there exists a unique homomorphism;
$\phi_i:SL_2(\bbC)\rightarrow G$ such that
\[
\hspace{-2pt}\phi_i\left(
\left(
\begin{array}{cc}
c&0\\
0&c^{-1}
\end{array}
\right)\right)=c^{\al^\vee_i},\,
\phi_i\left(
\left(
\begin{array}{cc}
1&t\\
0&1
\end{array}
\right)\right)=\exp(t e_i),\,
 \phi_i\left(
\left(
\begin{array}{cc}
1&0\\
t&1
\end{array}
\right)\right)=\exp(t f_i).
\]
where $c\in\bbC^\times$ and $t\in\bbC$.
Set $\al^\vee_i(c):=c^{\al^\vee_i}$,
$x_i(t):=\exp{(t e_i)}$, $y_i(t):=\exp{(t f_i)}$, 
$G_i:=\phi_i(SL_2(\bbC))$,
$T_i:=\phi_i(\{{\rm diag}(c,c^{-1})\vert 
c\in\bbC^{\vee}\})$ 
and 
$N_i:=N_{G_i}(T_i)$. Let
$T$ (resp. $N$) be the subgroup of $G$ 
with the Lie algebra $\tt$
(resp. generated by the $N_i$'s), 
which is called a {\it maximal torus} in $G$, and let
$B^{\pm}=U^{\pm}T$ be the Borel subgroup of $G$.
We have the isomorphism
$\phi:W\mapright{\sim}N/T$ defined by $\phi(s_i)=N_iT/T$.
An element $\ovl s_i:=x_i(-1)y_i(1)x_i(-1)
=\phi_i\left(
\left(
\begin{array}{cc}
0&\pm1\\
\mp1&0
\end{array}
\right)\right)$ is in 
$N_G(T)$, which is a representative of 
$s_i\in W=N_G(T)/T$. 

\subsection{Geometric crystals}
Let $W$ be the  Weyl group associated with $\ge$. 
Define $R(w)$ for $w\in W$ by
\[
 R(w):=\{(i_1,i_2,\cd,i_l)\in I^l|w=s_{i_1}s_{i_2}\cd s_{i_l}\},
\]
where $l$ is the length of $w$.
Then $R(w)$ is the set of reduced words of $w$.

Let $X$ be an ind-variety , 
{$\gamma_i:X\rightarrow \bbC$} and 
$\vep_i:X\longrightarrow \bbC$ ($i\in I$) 
rational functions on $X$, and
{$e_i:\bbC^\times \times X\longrightarrow X$}
$((c,x)\mapsto e^c_i(x))$ a
rational $\bbC^\times$-action.

For a word ${\bf i}=(i_1,\cd,i_l)\in R(w)$ 
$(w\in W)$, set 
$\al^{(j)}:=s_{i_l}\cd s_{i_{j+1}}(\al_{i_j})$ 
$(1\leq j\leq l)$ and 
\begin{eqnarray*}
e_{\bf i}:&T\times X\rightarrow &X\\
&(t,x)\mapsto &e_{\bf i}^t(x):=e_{i_1}^{\al^{(1)}(t)}
e_{i_2}^{\al^{(2)}(t)}\cd e_{i_l}^{\al^{(l)}(t)}(x).
\label{tx}
\end{eqnarray*}
\begin{df}
\label{def-gc}
A quadruple $(X,\{e_i\}_{i\in I},\{\gamma_i,\}_{i\in I},
\{\vep_i\}_{i\in I})$ is a 
$G$ (or $\ge$)-\\{\it geometric} {\it crystal} 
if
\begin{enumerate}
\item
$\{1\}\times X\subset dom(e_i)$ 
for any $i\in I$.
\item
$\gamma_j(e^c_i(x))=c^{a_{ij}}\gamma_j(x)$.
\item
{$e_{\bf i}=e_{\bf i'}$}
for any 
$w\in W$, ${\bf i}$.
${\bf i'}\in R(w)$.
\item
$\vep_i(e_i^c(x))=c^{-1}\vep_i(x)$.
\end{enumerate}
\end{df}
Note that the condition (iii) as above is 
equivalent to the following so-called 
{\it Verma relations}:
\[
 \begin{array}{lll}
&\hspace{-20pt}e^{c_1}_{i}e^{c_2}_{j}
=e^{c_2}_{j}e^{c_1}_{i}&
{\rm if }\,\,a_{ij}=a_{ji}=0,\\
&\hspace{-20pt} e^{c_1}_{i}e^{c_1c_2}_{j}e^{c_2}_{i}
=e^{c_2}_{j}e^{c_1c_2}_{i}e^{c_1}_{j}&
{\rm if }\,\,a_{ij}=a_{ji}=-1,\\
&\hspace{-20pt}
e^{c_1}_{i}e^{c^2_1c_2}_{j}e^{c_1c_2}_{i}e^{c_2}_{j}
=e^{c_2}_{j}e^{c_1c_2}_{i}e^{c^2_1c_2}_{j}e^{c_1}_{i}&
{\rm if }\,\,a_{ij}=-2,\,
a_{ji}=-1,\\
&\hspace{-20pt}
e^{c_1}_{i}e^{c^3_1c_2}_{j}e^{c^2_1c_2}_{i}
e^{c^3_1c^2_2}_{j}e^{c_1c_2}_{i}e^{c_2}_{j}
=e^{c_2}_{j}e^{c_1c_2}_{i}e^{c^3_1c^2_2}_{j}e^{c^2_1c_2}_{i}
e^{c^3_1c_2}_je^{c_1}_i&
{\rm if }\,\,a_{ij}=-3,\,
a_{ji}=-1,
\end{array}
\]
Note that the last formula is different from the one in 
\cite{BK}, \cite{N}, \cite{N2} which seems to be
incorrect. The formula here may be correct.

\subsection{Geometric crystal on Schubert cell}
\label{schubert}

Let $w\in W$ be a Weyl group element and take a 
reduced expression $w=s_{i_1}\cd s_{i_l}$. 
Let $X:=G/B$ be the flag
variety, which is an ind-variety 
and $X_w\subset X$ the
Schubert cell associated with $w$, which has 
a natural geometric crystal structure
(\cite{BK},\cite{N}).
For ${\bf i}:=(i_1,\cd,i_k)$, set 
\begin{equation}
B_{\bf i}^-
:=\{Y_{\bf i}(c_1,\cd,c_k)
:=Y_{i_1}(c_1)\cd Y_{i_l}(c_k)
\,\vert\, c_1\cd,c_k\in\bbC^\times\}\subset B^-,
\label{bw1}
\end{equation}
which has a geometric crystal structure(\cite{N})
isomorphic to $X_w$. 
The explicit forms of the action $e^c_i$, the rational 
function $\vep_i$  and $\gamma_i$ on 
$B_{\bf i}^-$ are given by
\begin{eqnarray}
&& e_i^c(Y_{i_1}(c_1)\cd Y_{i_l}(c_k))
=Y_{i_1}({\mathcal C}_1)\cd Y_{i_l}({\mathcal C}_k)),\nn \\
&&\text{where}\nn\\
&&{\mathcal C}_j:=
c_j\cdot \frac{\displaystyle \sum_{1\leq m\leq j,i_m=i}
 \frac{c}
{c_1^{a_{i_1,i}}\cd c_{m-1}^{a_{i_{m-1},i}}c_m}
+\sum_{j< m\leq k,i_m=i} \frac{1}
{c_1^{a_{i_1,i}}\cd c_{m-1}^{a_{i_{m-1},i}}c_m}}
{\displaystyle\sum_{1\leq m<j,i_m=i} 
 \frac{c}
{c_1^{a_{i_1,i}}\cd c_{m-1}^{a_{i_{m-1},i}}c_m}+
\mathop\sum_{j\leq m\leq k,i_m=i}  \frac{1}
{c_1^{a_{i_1,i}}\cd c_{m-1}^{a_{i_{m-1},i}}c_m}},
\label{eici}\\
&& \vep_i(Y_{i_1}(c_1)\cd Y_{i_l}(c_k))=
\sum_{1\leq m\leq k,i_m=i} \frac{1}
{c_1^{a_{i_1,i}}\cd c_{m-1}^{a_{i_{m-1},i}}c_m},
\label{vep-i}\\
&&\gamma_i(Y_{i_1}(c_1)\cd Y_{i_l}(c_k))
=c_1^{a_{i_1,i}}\cd c_k^{a_{i_k,i}}.
\label{gamma-i}
\end{eqnarray}

\subsection{Positive structure,\,\,
Ultra-discretizations \,\, and \,\,Tropicalizations}
\label{positive-str}
Let us recall the notions of 
positive structure, ultra-discretization and tropicalization.

The setting below is same as \cite{KNO}.
Let $T=(\bbC^\times)^l$ be an algebraic torus over $\bbC$ and 
$X^*(T):={\rm Hom}(T,\bbC^\times)\cong \ZZ^l$ 
(resp. $X_*(T):={\rm Hom}(\bbC^\times,T)\cong \ZZ^l$) 
be the lattice of characters
(resp. co-characters)
of $T$. 
Set $R:=\bbC(c)$ and define
$$
\begin{array}{cccc}
v:&R\setminus\{0\}&\longrightarrow &\ZZ\\
&f(c)&\mapsto
&{\rm deg}(f(c)),
\end{array}
$$
where $\rm deg$ is the degree of poles at $c=\ify$. 
Here note that for $f_1,f_2\in R\setminus\{0\}$, we have
\begin{equation}
v(f_1 f_2)=v(f_1)+v(f_2),\q
v\left(\frac{f_1}{f_2}\right)=v(f_1)-v(f_2)
\label{ff=f+f}
\end{equation}
A non-zero rational function on
an algebraic torus $T$ is called {\em positive} if
it is written as $g/h$ where
$g$ and $h$ are a positive linear combination of
characters of $T$.
\begin{df}
Let 
$f\cl T\rightarrow T'$ be 
a rational morphism between
two algebraic tori $T$ and 
$T'$.
We say that $f$ is {\em positive},
if $\chi\circ f$ is positive
for any character $\chi\cl T'\to \C$.
\end{df}
Denote by ${\rm Mor}^+(T,T')$ the set of 
positive rational morphisms from $T$ to $T'$.

\begin{lem}[\cite{BK}]
\label{TTT}
For any $f\in {\rm Mor}^+(T_1,T_2)$             
and $g\in {\rm Mor}^+(T_2,T_3)$, 
the composition $g\circ f$
is well-defined and belongs to ${\rm Mor}^+(T_1,T_3)$.
\end{lem}

By Lemma \ref{TTT}, we can define a category ${\mathcal T}_+$
whose objects are algebraic tori over $\bbC$ and arrows
are positive rational morphisms.

Let $f\cl T\rightarrow T'$ be a 
positive rational morphism
of algebraic tori $T$ and 
$T'$.
We define a map $\what f\cl X_*(T)\rightarrow X_*(T')$ by 
\[
\langle\chi,\what f(\xi)\rangle
=v(\chi\circ f\circ \xi),
\]
where $\chi\in X^*(T')$ and $\xi\in X_*(T)$.
\begin{lem}[\cite{BK}]
For any algebraic tori $T_1$, $T_2$, $T_3$, 
and positive rational morphisms 
$f\in {\rm Mor}^+(T_1,T_2)$, 
$g\in {\rm Mor}^+(T_2,T_3)$, we have
$\what{g\circ f}=\what g\circ\what f.$
\end{lem}
By this lemma, we obtain a functor 
\[
\begin{array}{cccc}
{\mathcal UD}:&{\mathcal T}_+&\longrightarrow &{{\hbox{\germ Set}}}\\
&T&\mapsto& X_*(T)\\
&(f:T\rightarrow T')&\mapsto& 
(\what f:X_*(T)\rightarrow X_*(T')))
\end{array}
\]


\begin{df}[\cite{BK}]
Let $\chi=(X,\{e_i\}_{i\in I},\{{\rm wt}_i\}_{i\in I},
\{\vep_i\}_{i\in I})$ be a 
geometric crystal, $T'$ an algebraic torus
and $\theta:T'\rightarrow X$ 
a birational isomorphism.
The isomorphism $\theta$ is called 
{\it positive structure} on
$\chi$ if it satisfies
\begin{enumerate}
\item for any $i\in I$ the rational functions
$\gamma_i\circ \theta:T'\rightarrow \bbC$ and 
$\vep_i\circ \theta:T'\rightarrow \bbC$ 
are positive.
\item
For any $i\in I$, the rational morphism 
$e_{i,\theta}:\bbC^\tm \tm T'\rightarrow T'$ defined by
$e_{i,\theta}(c,t)
:=\theta^{-1}\circ e_i^c\circ \theta(t)$
is positive.
\end{enumerate}
\end{df}
Let $\theta:T\rightarrow X$ be a positive structure on 
a geometric crystal $\chi=(X,\{e_i\}_{i\in I},$
$\{{\rm wt}_i\}_{i\in I},
\{\vep_i\}_{i\in I})$.
Applying the functor ${\mathcal UD}$ 
to positive rational morphisms
$e_{i,\theta}:\bbC^\tm \tm T'\rightarrow T'$ and
$\gamma\circ \theta:T'\ra T$
(the notations are
as above), we obtain
\begin{eqnarray*}
\til e_i&:=&{\mathcal UD}(e_{i,\theta}):
\ZZ\tm X_*(T) \rightarrow X_*(T)\\
{\rm wt}_i&:=&{\mathcal UD}(\gamma_i\circ\theta):
X_*(T')\rightarrow \bbZ,\\
\vep_i&:=&{\mathcal UD}(\vep_i\circ\theta):
X_*(T')\rightarrow \bbZ.
\end{eqnarray*}
Now, for given positive structure $\theta:T'\rightarrow X$
on a geometric crystal 
$\chi=(X,\{e_i\}_{i\in I},$
$\{{\rm wt}_i\}_{i\in I},
\{\vep_i\}_{i\in I})$, we associate 
the quadruple $(X_*(T'),\{\til e_i\}_{i\in I},
\{{\rm wt}_i\}_{i\in I},\{\vep_i\}_{i\in I})$
with a free pre-crystal structure (see \cite[2.2]{BK}) 
and denote it by ${\mathcal UD}_{\theta,T'}(\chi)$.
We have the following theorem:

\begin{thm}[\cite{BK}\cite{N}]
For any geometric crystal 
$\chi=(X,\{e_i\}_{i\in I},\{\gamma_i\}_{i\in I},$
$\{\vep_i\}_{i\in I})$ and positive structure
$\theta:T'\rightarrow X$, the associated pre-crystal 
${\mathcal UD}_{\theta,T'}(\chi)=
(X_*(T'),\{e_i\}_{i\in I},\{{\rm wt}_i\}_{i\in I},
\{\vep_i\}_{i\in I})$ 
is a crystal {\rm (see \cite[2.2]{BK})}
\end{thm}

Now, let ${\mathcal GC}^+$ be a category whose 
object is a triplet
$(\chi,T',\theta)$ where 
$\chi=(X,\{e_i\},\{\gamma_i\},\{\vep_i\})$ 
is a geometric crystal and $\theta:T'\rightarrow X$ 
is a positive structure on $\chi$, and morphism
$f:(\chi_1,T'_1,\theta_1)\longrightarrow 
(\chi_2,T'_2,\theta_2)$ is given by a morphism 
$\vp:X_1\longrightarrow X_2$  
($\chi_i=(X_i,\cd)$) such that 
\[
f:=\theta_2^{-1}\circ\vp\circ\theta_1:T'_1\longrightarrow T'_2,
\]
is a positive rational morphism. Let ${\mathcal CR}$
be a category of crystals. 
Then by the theorem above, we have
\begin{cor}
\label{cor-posi}
$\mathcal UD_{\theta,T'}$ as above defines a functor
\begin{eqnarray*}
 {\mathcal UD}&:&{\mathcal GC}^+\longrightarrow {\mathcal CR},\\
&&(\chi,T',\theta)\mapsto X_*(T'),\\
&&(f:(\chi_1,T'_1,\theta_1)\rightarrow 
(\chi_2,T'_2,\theta_2))\mapsto
(\what f:X_*(T'_1)\rightarrow X_*(T'_2)).
\end{eqnarray*}

\end{cor}
We call the functor $\mathcal UD$
{\it ``ultra-discretization''} as \cite{N},\cite{N2}
instead of ``tropicalization'' as in \cite{BK}.
And 
for a crystal $B$, if there
exists a geometric crystal $\chi$ and a positive 
structure $\theta:T'\rightarrow X$ on $\chi$ such that 
${\mathcal UD}(\chi,T',\theta)\cong B$ as crystals, 
we call an object $(\chi,T',\theta)$ in ${\mathcal GC}^+$
a {\it tropicalization} of $B$, where 
it is not known that this correspondence is a functor.

\renewcommand{\thesection}{\arabic{section}}
\section{Limit of perfect crystals}
\label{limit}
\setcounter{equation}{0}
\renewcommand{\theequation}{\thesection.\arabic{equation}}
We review limit of perfect crystals following \cite{KKM}.
(See also \cite{KMN1},\cite{KMN2}).

\subsection{Crystals}

First we review the theory of crystals,
which is the notion obtained by
abstracting the combinatorial 
properties of crystal bases.
\begin{df}
A {\it crystal} $B$ is a set endowed with the following maps:
\begin{eqnarray*}
&& {\rm wt}:B\lar P,\\
&&\vep_i:B\lar\ZZ\sqcup\{-\infty\},\q
  \vp_i:B\lar\ZZ\sqcup\{-\infty\} \q{\hbox{for}}\q i\in I,\\
&&\eit:B\sqcup\{0\}\lar B\sqcup\{0\},
\q\fit:B\sqcup\{0\}\lar B\sqcup\{0\}\q{\hbox{for}}\q i\in I,\\
&&\eit(0)=\fit(0)=0.
\end{eqnarray*}
those maps satisfy the following axioms: for
 all $b,b_1,b_2 \in B$, we have
\begin{eqnarray*}
&&\vp_i(b)=\vep_i(b)+\lan \al^\vee_i,{\rm wt}
(b)\ran,\\
&&\wt(\eit b)=\wt(b)+\al_i{\hbox{ if  }}\eit b\in B,\\
&&\wt(\fit b)=\wt(b)-\al_i{\hbox{ if  }}\fit b\in B,\\
&&\eit b_2=b_1 \Longleftrightarrow \fit b_1=b_2\,\,(\,b_1,b_2 \in B),\\
&&\vep_i(b)=-\ify
   \Longrightarrow \eit b=\fit b=0.
\end{eqnarray*}
\end{df}
The following tensor product structure 
is one of the most crucial properties of crystals.
\begin{thm}
\label{tensor}
Let $B_1$ and $B_2$ be crystals.
Set
$B_1\ot B_2:=
\{b_1\otimes b_2;\;b_j\in B_j\;(j=1,2)\}$. Then we have 
\begin{enumerate}
\item $B_1\ot B_2$ is a crystal.
\item
For $b_1\in B_1$ and $b_2\in B_2$, we have
$$
\tilde f_i(b_1\otimes b_2)=
\left\{\begin{array}{ll}\tilde f_ib_1\otimes b_2&
{\rm if}\;\varphi_i(b_1)>\vep_i(b_2),\\
b_1\otimes\tilde f_ib_2&{\rm if}\;
\varphi_i(b_1)\leq\vep_i(b_2).
\end{array}\right.
$$
$$
\tilde e_i(b_1\otimes b_2)=\left\{\begin{array}{ll}
b_1\otimes \tilde e_ib_2&
{\rm if}\;\varphi_i(b_1)<\vep_i(b_2),\\
\tilde e_ib_1\otimes b_2
&{\rm if}\;\varphi_i(b_1)\geq\vep_i(b_2),
\end{array}\right.
$$
\end{enumerate}
\end{thm}

\begin{df}
Let $B_1$ and $B_2$ be crystals. A {\it strict morphism} of crystals
$\psi:B_1\lar B_2$ is a map
$\psi:B_1\sqcup\{0\} \lar B_2\sqcup\{0\}$
satisfying: 
$\psi(0)=0$, $\psi(B_1)\subset B_2$,
$\psi$ commutes with all $\eit$ and $\fit$
and
\[
\hspace{-30pt}\wt(\psi(b))=\wt(b),\q \vep_i(\psi(b))=\vep_i(b),\q
  \vp_i(\psi(b))=\vp_i(b)
\text{ for any }b\in B_1.
\]
In particular, 
a bijective strict morphism is called an 
{\it isomorphism of crystals}. 
\end{df}

\begin{ex}
\label{ex-tlm}
If $(L,B)$ is a crystal base, then $B$ is a crystal.
Hence, for the crystal base $(L(\ify),B(\ify))$
of the nilpotent subalgebra $\uqm$ of 
the quantum algebra $\uq$, 
$B(\ify)$ is a crystal. 
\end{ex}
\begin{ex}
\label{tlm}
For $\lm\in P$, set $T_\lm:=\{t_\lm\}$. We define a crystal
structure on $T_\lm$ by 
\[
 \eit(t_\lm)=\fit(t_\lm)=0,\q\vep_i(t_\lm)=
\vp_i(t_\lm)=-\ify,\q \wt(t_\lm)=\lm.
\]
\end{ex}
\begin{df}
For a crystal $B$, a colored oriented graph
structure is associated with $B$ by 
\[
 b_1\mapright{i}b_2\Longleftrightarrow 
\fit b_1=b_2.
\]
We call this graph a {\it crystal graph}
of $B$.
\end{df}

\subsection{Affine weights}
\label{aff-wt}

Let $\ge$ be an affine Lie algebra. 
The sets $\mathfrak t$, 
$\{\al_i\}_{i\in I}$ 
and $\{\al^\vee_i\}_{i\in I}$ be as in \ref{KM}. 
We take ${\rm dim}\mathfrak t=\sharp I+1$.
Let $\del\in Q_+$ be the unique element 
satisfying $\{\lm\in Q|\lan \al^\vee_i,\lm\ran=0
\text{ for any }i\in I\}=\bbZ\del$
and ${\bf c}\in \ge$ be the canonical central element
satisfying $\{h\in Q^\vee|\lan h,\al_i\ran=0
\text{ for any }i\in I\}=\bbZ c$.
We write (\cite[6.1]{Kac})
\[
{\bf c}=\sum_i a_i^\vee \al^\vee_i,\qq
\del=\sum_i a_i\al_i.
\]
Let $(\q,\q)$ be the non-degenerate
$W$-invariant symmetric bilinear form on $\mathfrak t^*$
normalized by $(\del,\lm)=\lan {\bf c},\lm\ran$
for $\lm\in\frak t^*$.
Let us set $\tt^*_{\rm cl}:=\tt^*/\bbC\del$ and let
${\rm cl}:\tt^*\longrightarrow \tt^*_{\rm cl}$
be the canonical projection. 
Here we have 
$\tt^*_{\rm cl}\cong \oplus_i(\bbC \al^\vee_i)^*$.
Set $\tt^*_0:=\{\lm\in\tt^*|\lan {\bf c},\lm\ran=0\}$,
$(\tt^*_{\rm cl})_0:={\rm cl}(\tt^*_0)$. 
Since $(\del,\del)=0$, we have a positive-definite
symmetric form on $\tt^*_{\rm cl}$ 
induced by the one on 
$\tt^*$. 
Let $\Lm_i\in \tt^*_{\rm cl}$ $(i\in I)$ be a classical 
weight such that $\lan \al^\vee_i,\Lm_j\ran=\del_{i,j}$, which 
is called a fundamental weight.
We choose 
$P$ so that $P_{\rm cl}:={\rm cl}(P)$ coincides with 
$\oplus_{i\in I}\bbZ\Lm_i$ and 
we call $P_{\rm cl}$ a 
{\it classical weight lattice}.

\subsection{Definitions of perfect crystal and its limit}
\label{def-perfect}

Let $\ge$ be an affine Lie algebra, $P_{cl}$ be
a classical weight lattice as above and set 
$(P_{cl})^+_l:=\{\lm\in P_{cl}|
\lan c,\lm\ran=l,\,\,\lan \al^\vee_i,\lm\ran\geq0\}$ 
$(l\in\ZZ_{>0})$.
\begin{df}
\label{perfect-def}
A crystal $B$ is a {\it perfect} of level $l$ if 
\begin{enumerate}
\item
$B\ot B$ is connected as a crystal graph.
\item
There exists $\lm_0\in P_{\rm cl}$ such that 
\[
 \wt(B)\subset \lm_0+\sum_{i\ne0}\ZZ_{\leq0}
{\rm cl}(\al_i),\qq
\sharp B_{\lm_0}=1
\]
\item There exists a finite-dimensional 
$U'_q(\ge)$-module $V$ with a
crystal pseudo-base $B_{ps}$ 
such that $B\cong B_{ps}/{\pm1}$
\item
The maps 
$\vep,\vp:B^{min}:=\{b\in B|\lan c,\vep(b)\ran=l\}
\mapright{}(P_{\rm cl}^+)_l$ are bijective, where 
$\vep(b):=\sum_i\vep_i(b)\Lm_i$ and 
$\vp(b):=\sum_i\vp_i(b)\Lm_i$.
\end{enumerate}
\end{df}

Let $\{B_l\}_{l\geq1}$ be a family of 
perfect crystals of level $l$ and set 
$J:=\{(l,b)|l>0,\,b\in B^{min}_l\}$.
\begin{df}
\label{def-limit}
A crystal $B_\ify$ with an element $b_\ify$ is called a
{\it limit of $\{B_l\}_{l\geq1}$}
if 
\begin{enumerate}
\item
$\wt(b_\ify)=\vep(b_\ify)=\vp(b_\ify)=0$.
\item
For any $(l,b)\in J$, there exists an
embedding of 
crystals:
\begin{eqnarray*}
 f_{(l,b)}:&
T_{\vep(b)}\ot B_l\ot T_{-\vp(b)}\hookrightarrow
B_\ify\\
&t_{\vep(b)}\ot b\ot t_{-\vp(b)}\mapsto b_\ify
\end{eqnarray*}
\item
$B_\ify=\bigcup_{(l,b)\in J} {\rm Im}f_{(l,b)}$.
\end{enumerate}
\end{df}
\noindent
As for the crystal $T_\lm$, see Example \ref{tlm}.
If a limit exists for a family $\{B_l\}$, 
we say that $\{B_l\}$
is a {\it coherent family} of perfect crystals.

The following is one of the most 
important properties of limit of perfect crystals.
\begin{pro}
Let $B(\ify)$ be the crystal as in 
Example \ref{ex-tlm}. Then we have
the following isomorphism of crystals:
\[
B(\ify)\ot B_\ify\mapright{\sim}B(\ify).
\]
\end{pro}

\renewcommand{\thesection}{\arabic{section}}
\section{Fundamental Representations}
\setcounter{equation}{0}
\renewcommand{\theequation}{\thesection.\arabic{equation}}

\subsection{Fundamental representation 
$W(\varpi_1)$}
\label{fundamental}

Let $c=\sum_{i}a_i^\vee \al^\vee_i$ be the canonical
central element in an affine Lie algebra $\ge$
(see \cite[6.1]{Kac}), 
$\{\Lm_i|i\in I\}$ the set of fundamental 
weight as in the previous section
and $\varpi_1:=\Lm_1-a^\vee_1\Lm_0$ the
(level 0)fundamental weight.

Let $V(\varpi_1)$ be the extremal weight module
of $\uq$
associated with $\varpi_1$ (\cite{K0}) and 
$W(\varpi_1)\cong 
V(\varpi_1)/(z_1-1)V(\varpi_1)$ the 
fundamental representation of $\uqp$
where $z_1$ is a 
$\uqp$-linear automorhism on $V(\varpi_1)$
 (see \cite[Sect 5.]{K0}). 

By \cite[Theorem 5.17]{K0}, $W(\varpi_1)$ is an
finite-dimensional irreducible integrable 
$\uqp$-module and has a global basis
with a simple crystal. Thus, we can consider 
the specialization $q=1$ and obtain the 
finite-dimensional $\ge$-module $W(\varpi_1)$, 
which we call a fundamental representation
of $\ge$ and use the same notation as above.

We shall present the explicit form of 
$W(\varpi_1)$ for $\ge=\TY(G,1,2)$.
\subsection{$W(\varpi_1)$ for $\TY(G,1,2)$}
The Cartan matrix $A=(a_{i,j})_{i,j=0,1,2}$ of type 
$\TY(G,1,2)$ is as follows:
\[
 A=\begin{pmatrix}2&-1&0\\
-1&2&-1\\0&-3&2
\end{pmatrix}.
\]
Then the simple roots are 
\[
 \al_0=2\Lm_0-\Lm_1+\del,\q
\al_1=-\Lm_0+2\Lm_1-3\Lm_2,\q
\al_2=-\Lm_1+2\Lm_2, 
\]
and the Dynkin diagram is:
\[\SelectTips{cm}{}
\xymatrix{
*{\bigcirc}<3pt> \ar@{-}[r]_<{0} 
& *{\bigcirc}<3pt> \ar@3{->}[r]_<{1}
& *{\bigcirc}<6pt>\ar@{}_<{\,\,\,\,\,\,2}
}
\]

The $\ge$-module $W(\varpi_1)$ is a 15 dimensional 
module with the basis,
\[
 \{\fsquare(5mm,i),\fsquare(5mm,\ovl i),
\emptyset, \fsquare(5mm,0_1),
\fsquare(5mm,0_2) \,\,\vert \,\,i=1,\cd,6\}.
\]
The explicit form of $W(\varpi_1)$ is given 
in \cite{Y}, which slightly differs from
our description  below.
\begin{eqnarray*}
&&\hspace{-30pt}{\rm wt}(\fsquare(5mm,1))=\Lm_1-2\Lm_0,\,\,
{\rm wt}(\fsquare(5mm,2))=-\Lm_0-\Lm_1+3\Lm_2,\,\,
{\rm wt}(\fsquare(5mm,3))=-\Lm_0+\Lm_2,\\
&&\hspace{-30pt}{\rm wt}(\fsquare(5mm,4))=-\Lm_0+\Lm_1-\Lm_2,\,\,
{\rm wt}(\fsquare(5mm,5))=-\Lm_1+2\Lm_2,\,\,
{\rm wt}(\fsquare(5mm,6))=-\Lm_0+2\Lm_1-3\Lm_2,\\
&&\hspace{-30pt}{\rm wt}(\fsquare(5mm,\ovl i))=
-{\rm wt}(\fsquare(5mm,i))\,\,(i=1,\cd,6),\,\,
{\rm wt}(\fsquare(5mm,0_1))=
{\rm wt}(\fsquare(5mm,0_2))=
{\rm wt}(\emptyset)=0.
\end{eqnarray*}
\def\bv#1{\fsquare(5mm,#1)}
The actions of $e_i$ and $f_i$ on these basis vectors
are given as follows:
\begin{eqnarray*}
&&\hspace{-30pt}
f_0\left(\bv{0_2},\bv{\ovl 6},\bv{\ovl 4},\bv{\ovl 3},\bv{\ovl 2},
\bv{\ovl 1},\emptyset\right)
=\left(\bv1,\bv2,\bv3,\bv4,\bv6,\emptyset,2\bv1\right),\\
&&\hspace{-30pt}e_0\left(\bv1,\bv2,\bv3,\bv4,\bv6,\bv{0_2},\emptyset\right)
=\left(\emptyset,\bv{\ovl 6},\bv{\ovl 4},\bv{\ovl 3},
\bv{\ovl 2},\bv{\ovl 1},2\bv{\ovl 1}\right),\\
&&\hspace{-30pt}f_1\left(\bv1,\bv4,\bv6,\bv{0_1},\bv{0_2},
\bv{\ovl 5},\bv{\ovl 2},\emptyset\right)=
\left(\bv2,\bv5,\bv{0_2},3\bv{\ovl 6},2\bv{\ovl 6},
\bv{\ovl 4},\bv{\ovl 1},\bv{\ovl 6}\right),\\
&&\hspace{-30pt}e_1\left(\bv2,\bv5,\bv{0_1},\bv{0_2},\bv{\ovl 6},
\bv{\ovl 4},\bv{\ovl 1},\emptyset\right)=
\left(\bv1,\bv4,3\bv6,2\bv6,\bv{0_2},\bv{\ovl 5},
\bv{\ovl 2},\bv6\right),\\
&&\hspace{-30pt}f_2\left(\bv2,\bv3,\bv4,\bv5,\bv{0_1},\bv{0_2},
\bv{\ovl 6},\bv{\ovl 4},\bv{\ovl 3}\right)\\
&&\qq\qq\qq\qq\qq\qq=
\left(\bv3,2\bv4,3\bv6,\bv{0_1},2\bv{\ovl 5},\bv{\ovl 5},
\bv{\ovl 4},2\bv{\ovl 3},3\bv{\ovl 2}\right),\\
&&\hspace{-30pt}e_2\left(\bv3,\bv4,\bv6,\bv{0_1},\bv{0_2},\bv{\ovl 5},
\bv{\ovl 4},\bv{\ovl 3},\bv{\ovl 2}\right)\\
&&\qq\qq\qq\qq\qq\qq=
\left(3\bv2,2\bv3,\bv4,2\bv5,\bv5,\bv{0_1},
3\bv{\ovl 6},2\bv{\ovl 4},\bv{\ovl 3}\right),
\end{eqnarray*}
where we only give non-trivial actions 
and the other actions are trivial.
We can easily check that these define 
the module $W(\varpi_1)$ by direct calculations.

\renewcommand{\thesection}{\arabic{section}}
\section{Affine Geometric Crystal $\cV_1(\TY(G,1,2))$}
\setcounter{equation}{0}
\renewcommand{\theequation}{\thesection.\arabic{equation}}

We shall construct the 
affine geometric crystal $\cV(\TY(G,1,2))$ in $W(\varpi_1)$
explicitly.

For $\xi\in (\frt^*_{\rm cl})_0$, let $t(\xi)$ be the 
shift as in \cite[Sect 4]{K0}.
Then we have 
\begin{eqnarray*}
&& t(\wtil\varpi_1)=s_0s_1s_2s_1s_2s_1=:w_1,\\
&& t(\text{wt}(\bv{\ovl 2}))=s_2s_1s_2s_1s_0s_1=:w_2,
\end{eqnarray*}
Associated with these Weyl group elements $w_1$ and $w_2$,
we define algebraic varieties $\cV_1=\cV_1(\TY(G,1,2))$ and 
$\cV_2=\cV_2(\TY(G,1,2))\subset W(\varpi_1)$ respectively:
\begin{eqnarray*}
&&\hspace{-30pt}\cV_1:=\{v_1(x)
:=Y_0(x_0)Y_1(x_1)Y_2(x_2)Y_1(x_3)Y_2(x_4)Y_1(x_5)
\bv1\,\,\vert\,\,x_i\in\bbC^\times,(0\leq i\leq 5)\},\\
&&\hspace{-30pt}\cV_2:=\{v_2(y):=
Y_2(y_2)Y_1(y_1)Y_2(y_4)Y_1(y_3)Y_0(y_0)Y_1(y_5)
\bv{\ovl 2}\,\,\vert\,\,y_i\in\bbC^\times,(0\leq i\leq 5)\}.
\end{eqnarray*}
Due to the explicit forms of $f_i$'s on $W(\varpi_1)$
as above, we have $f_0^3=0$, $f_1^3=0$ and $f_2^4=0$ 
and then 
\begin{equation}
Y_i(c)=(1+\frac{f_i}{c}+\frac{f_i^2}{2c^2})\al_i^\vee(c)
\,\,(i=0,1),\q
Y_2(c)=(1+\frac{f_2}{c}+\frac{f_2^2}{2c^2}
+\frac{f_2^3}{6c^3})\al_2^\vee(c).
\end{equation}
Thus, we can get explicit forms of $v_1(x)\in\cV_1$ 
and $v_2(y)\in\cV_2$. 
Set 
\begin{eqnarray*}
&&v_1(x)=\sum_{1\leq i\leq 6}\left(X_i\bv{i}+X_{\ovl i}
\bv{\ovl i}\right)+X_{0_1}\bv{0_1}+X_{0_2}\bv{0_2}+X_\emptyset
\emptyset,\\
&&v_2(y)=\sum_{1\leq i\leq 6}\left(Y_i\bv{i}+Y_{\ovl i}
\bv{\ovl i}\right)+Y_{0_1}\bv{0_1}+Y_{0_2}\bv{0_2}+Y_\emptyset
\emptyset.
\end{eqnarray*}
Then by direct calculations, we have 
\begin{lem}
\label{XY}
The rational function $X_1,X_2,\cd,$ and $Y_1,Y_2,\cd$
are given as:
\begin{eqnarray*}
&&\hspace{-20pt}
X_1=1 + \frac{{x_3}}{{x_0}} + \frac{{x_1}\,{{x_3}}^2}
   {{x_0}\,{{x_2}}^3} + \frac{3\,{x_1}\,{x_3}\,{x_4}}
   {{x_0}\,{{x_2}}^2} + \frac{3\,{x_1}\,{{x_4}}^2}
   {{x_0}\,{x_2}} + \frac{{x_1}\,{{x_4}}^3}{{x_0}\,{x_3}} + 
  \left( \frac{{x_1}}{{x_0}} 
+ \frac{{x_1}\,{x_3}}{{{x_0}}^2} \right)\,{x_5},\\
&&\hspace{-20pt}X_2=
\frac{{{x_2}}^3}{{{x_1}}^2} + \frac{{{x_3}}^2}{{{x_2}}^3} + 
  \frac{3\,{x_3}\,{x_4}}{{{x_2}}^2} 
+ \frac{3\,{{x_4}}^2}{{x_2}} + 
  \frac{{{x_4}}^3}{{x_3}} + {x_5} + 
  \frac{{x_3}\,{x_5}}{{x_0}}  \\
&&\hspace{-20pt}\qq\q
+\frac{{x_0}\,{x_3}\,\left( 2\,{x_3} + 3\,{x_2}\,{x_4} \right)+ 
  {{x_2}}^3\,\left( {{x_4}}^3 + {x_3}\,{x_5} \right) }{{x_0}\,
 {x_1}\,{x_3}},\\
&&\hspace{-20pt}X_3=
\frac{{{x_2}}^2}{{x_1}} + \frac{{x_3}}{{x_2}} + {x_4} + 
  \frac{{x_2}\,{{x_4}}^2}{{x_0}} + 
  \frac{{{x_2}}^2\,{{x_4}}^3}{{x_0}\,{x_3}} + 
  \frac{{{x_2}}^2\,{x_5}}{{x_0}},\\
&&\hspace{-20pt}X_4=
{x_2} + \frac{{x_1}\,{x_3}\,{x_4}}{{x_0}\,{x_2}} + 
  \frac{2\,{x_1}\,{{x_4}}^2}{{x_0}} + 
  \frac{{x_1}\,{x_2}\,{{x_4}}^3}{{x_0}\,{x_3}} + 
  \frac{{x_1}\,{x_2}\,{x_5}}{{x_0}},\\
&&\hspace{-20pt}X_5=
\left( \frac{{{x_2}}^2}{{x_1}} + \frac{{x_3}}{{x_2}} \right) \,
   {x_4} + 2\,{{x_4}}^2 + \frac{{x_2}\,{{x_4}}^3}{{x_3}} + 
  {x_2}\,{x_5},\\
&&\hspace{-20pt}X_6=
{x_1} + \frac{{{x_1}}^2\,{{x_3}}^2}{{x_0}\,{{x_2}}^3} + 
  \frac{3\,{{x_1}}^2\,{x_3}\,{x_4}}{{x_0}\,{{x_2}}^2} + 
  \frac{3\,{{x_1}}^2\,{{x_4}}^2}{{x_0}\,{x_2}} + 
  \frac{{{x_1}}^2\,{{x_4}}^3}{{x_0}\,{x_3}} + 
  \frac{{{x_1}}^2\,{x_5}}{{x_0}},\\
&&\hspace{-20pt}X_{0_1}=x_2x_4,\q
X_{0_2}=
{x_3} + \frac{{x_1}\,{{x_3}}^2}{{{x_2}}^3} + 
  \frac{3\,{x_1}\,{x_3}\,{x_4}}{{{x_2}}^2} + 
  \frac{3\,{x_1}\,{{x_4}}^2}{{x_2}} + 
  \frac{{x_1}\,{{x_4}}^3}{{x_3}} + {x_1}\,{x_5},\\
&&\hspace{-20pt}X_{\ovl 6}=
{x_0}\,\left( \frac{{{x_2}}^3}{{{x_1}}^2} + 
    \frac{{{x_3}}^2}{{{x_2}}^3} 
   + \frac{3\,{x_3}\,{x_4}}{{{x_2}}^2} + 
    \frac{3\,{{x_4}}^2}{{x_2}} + \frac{{{x_4}}^3}{{x_3}} + 
    \frac{2\,{x_3} + 3\,{x_2}\,{x_4}}{{x_1}} + {x_5} \right),\\
&&\hspace{-20pt}X_{\ovl 5}=
{x_1}\,\left( \frac{{x_3}}{{x_2}} + {x_4} \right),\q
X_{\ovl 4}={x_0}\,\left( \frac{{{x_2}}^2}{{x_1}} 
+ \frac{{x_3}}{{x_2}} + {x_4} \right),\q
X_{\ovl 3}=x_0x_2,\\
&&\hspace{-20pt}X_{\ovl 2}=x_0x_1,\q
X_{\ovl 1}=x_0^2,\q
X_\emptyset=x_0,\\
&&\hspace{-20pt}Y_1=y_1y_3,
\q Y_2=\frac{{{y_2}}^3\,\left( {y_1}\,{y_3} 
+ {{y_4}}^3 \right) }{{y_1}},\q
Y_3={{y_2}}^2\,{y_3} + {y_2}\,{{y_4}}^2 + 
  \frac{{{y_2}}^2\,{{y_4}}^3}{{y_1}},
\\
&&\hspace{-20pt}Y_4={y_2}\,{y_3} + \frac{{y_1}\,{y_4}}{{y_2}} + 2\,{{y_4}}^2 + 
  \frac{{y_2}\,{{y_4}}^3}{{y_1}},\q
Y_5=y_2^2y_4,\\
&&\hspace{-20pt}Y_6={y_3} + \frac{3\,{y_1}\,{y_4}}{{{y_2}}^2} + 
  \frac{3\,{{y_4}}^2}{{y_2}} + \frac{{{y_4}}^3}{{y_1}} + 
  \frac{{{y_1}}^2\,\left( {{y_4}}^3 + {{y_3}}^2\,{y_5} \right) }
   {{{y_2}}^3\,{{y_4}}^3},\\
&&\hspace{-20pt}Y_{0_1}=y_2y_4,\q
Y_{0_2}={y_1} + \left( {y_3} 
+ \frac{{y_1}\,{{y_3}}^2}{{{y_4}}^3} \right) \,
   {y_5},\\
&&\hspace{-20pt}Y_{\ovl 6}=
{{y_2}}^3 + {y_0}\,\left( \frac{{{y_2}}^3}{{y_1}} + 
     \frac{{{y_2}}^3\,{{y_4}}^3}{{{y_1}}^2\,{y_3}} \right)  + 
  \left( \frac{2\,{{y_2}}^3\,{y_3}}{{y_1}} + 
     \frac{{{y_2}}^3\,{{y_3}}^2}{{{y_4}}^3} + 
     \frac{{{y_2}}^3\,{{y_4}}^3}{{{y_1}}^2} \right) \,{y_5},\\
&&\hspace{-20pt}Y_{\ovl 5}=
\frac{{y_1}}{{y_2}} + {y_4} + 
  \left( {y_3}\,\left( \frac{1}{{y_2}} + 
        \frac{{y_1}}{{{y_2}}^2\,{y_4}} \right)  + 
     \frac{{y_1}\,{{y_3}}^2}{{y_2}\,{{y_4}}^3} \right) \,{y_5},
\\
&&\hspace{-20pt}Y_{\ovl 4}=
{{y_2}}^2 
+ {y_0}\,\left( \frac{{{y_2}}^2}{{y_1}} + 
     \frac{{y_2}\,{{y_4}}^2}{{y_1y_3}} + 
     \frac{{{y_2}}^2\,{{y_4}}^3}{{{y_1}}^2\,{y_3}} \right) \\
&& + 
  \left( {y_3}\,\left( \frac{2\,{{y_2}}^2}{{y_1}} + 
        \frac{{y_2}}{{y_4}} \right)  + 
     \frac{{{y_2}}^2\,{{y_3}}^2}{{{y_4}}^3} + 
     \frac{{y_2}\,{{y_4}}^2}{{y_1}} + 
\frac{{{y_2}}^2\,{{y_4}}^3}{{{y_1}}^2} \right) \,{y_5}, \\
&&\hspace{-20pt} Y_{\ovl 3}=
{y_2} + {y_0}\,\left( \frac{{y_2}}{{y_1}} + 
   \frac{{y_4}}{{y_2y_3}} + \frac{2\,{{y_4}}^2}{{y_1y_3}} + 
   \frac{{y_2}\,{{y_4}}^3}{{{y_1}}^2\,{y_3}} \right) \\ 
&&\hspace{-20pt}\qq\q 
+\left( {y_3}\,\left( \frac{2\,{y_2}}{{y_1}} 
+ \frac{2}{{y_4}} \right)  + 
\frac{{y_2}\,{{y_3}}^2}{{{y_4}}^3} + \frac{{y_4}}{{y_2}} + 
\frac{2\,{{y_4}}^2}{{y_1}} + \frac{{y_2}\,{{y_4}}^3}{{{y_1}}^2}
     \right) \,{y_5},
\end{eqnarray*}
\begin{eqnarray*}
&&\hspace{-20pt}Y_{\ovl 2}=
1 + {y_0}\,\left( \frac{1}{{y_1}} 
+ \frac{{y_1}}{{{y_2}}^3\,{y_3}} + 
     \frac{3\,{y_4}}{{{y_2}}^2\,{y_3}} + 
     \frac{3\,{{y_4}}^2}{{y_1}\,{y_2y_3}} + 
     \frac{{{y_4}}^3}{{{y_1}}^2\,{y_3}} \right) \\
&&\hspace{-20pt}\qq\q + 
  \left( \frac{{y_1}}{{{y_2}}^3} + 
{y_3}\,\left( \frac{2}{{y_1}} + \frac{3}{{y_2}\,{y_4}} \right)  + 
\frac{{y_0}\,{y_1}\,{y_3}}{{{y_2}}^3\,{{y_4}}^3} + 
\frac{{{y_3}}^2}{{{y_4}}^3} + \frac{3\,{y_4}}{{{y_2}}^2} + 
\frac{3\,{{y_4}}^2}{{y_1}\,{y_2}} + \frac{{{y_4}}^3}{{{y_1}}^2}
     \right) \,{y_5},\\
&&\hspace{-20pt}Y_{\ovl 1}=
\frac{{{y_0}}^2}{{y_1}\,{y_3}} + {y_5} + 
  {y_0}\,\left( \frac{1}{{y_3}} + \frac{{y_5}}{{y_1}} + 
     \frac{{y_3}\,{y_5}}{{{y_4}}^3} \right),\q
Y_\emptyset=y_0.
\end{eqnarray*}
\end{lem}
Now we solve the equation 
\begin{equation}
 v_2(y)=a(x)v_1(x),
\label{eq}
\end{equation}
where $a(x)$ is a rational function in $x=(x_0,\cd,x_6)$.
Though this equation is over-determined, 
we can solve it and obtain the explicit form of the unique
solution as follows:
\begin{pro}
Set
\begin{eqnarray*}
&&M:=\frac{{x_3}\,{{x_4}}^2}{{x_0}\,{x_2}} + \frac{3\,{{x_4}}^3}{{x_0}} + 
  \frac{{x_3}\,{x_5}}{{x_2}\,{x_4}} + 
  {x_2}\,\left( \frac{3\,{{x_4}}^4}{{x_0}\,{x_3}} + 
     \frac{3\,{x_4}\,{x_5}}{{x_0}} \right)\\
&&\qq  + 
  {{x_2}}^2\,\left( \frac{{{x_4}}^2}{{x_0}\,{x_1}} + 
     \frac{{{x_4}}^2}{{x_1}\,{x_3}} 
 + \frac{{{x_4}}^5}{{x_0}\,{{x_3}}^2} + 
     \frac{{x_5}}{{x_1}\,{x_4}} + 
     \frac{2\,{{x_4}}^2\,{x_5}}{{x_0}\,{x_3}} + 
     \frac{{{x_5}}^2}{{x_0}\,{x_4}} \right),\\
&&N:=
\frac{3\,{x_1}\,{x_3}}{{{x_2}}^3} + 
  \frac{{x_2}\,{x_3}}{{x_1}\,{{x_4}}^2} + 
  \frac{2\,{{x_3}}^2}{{{x_2}}^2\,{{x_4}}^2} + 
  \frac{{x_1}\,{{x_3}}^3}{{{x_2}}^5\,{{x_4}}^2} + 
  \frac{3\,{x_3}}{{x_2}\,{x_4}} + 
  \frac{3\,{x_1}\,{{x_3}}^2}{{{x_2}}^4\,{x_4}}\\ 
&&\qq\q 
+\frac{{x_1}\,{x_4}}{{{x_2}}^2} + \frac{{x_2}\,{x_4}}{{x_0}} + 
  \frac{{x_1}\,{x_3}\,{x_5}}{{{x_2}}^2\,{{x_4}}^2} + 
  \frac{{x_2}\,{x_3}\,{x_5}}{{x_0}\,{{x_4}}^2} + 
  \frac{{x_1}\,{{x_3}}^2\,{x_5}}{{x_0}\,{{x_2}}^2\,{{x_4}}^2}.
\end{eqnarray*}
Then we have the rational function $a(x)$ and 
the unique solution of (\ref{eq}):
\begin{equation}
\begin{split}
&a(x)=\frac{M}{(x_2x_4)^2},\q
y_2=\frac{{x_2}}{{x_1}} + \frac{{x_3}}{{{x_2}}^2} + 
  \frac{2\,{x_4}}{{x_2}} + \frac{{{x_4}}^2}{{x_3}} + 
  \frac{{x_5}}{{x_4}},\q
y_4=\frac{M}{y_2x_2x_4},\\
&y_0=a(x)x_0,\,\,
y_1=\frac{y_2^3(a(x)X_1+y_4^3)}{a(x)X_2},\q
y_3=\frac{a(x)X_1}{y_1},\q
y_5=\frac{\frac{x_5MN}{x_1x_2x_3x_4}}
{a(x)X_2\left(y_3+\frac{y_1y_3^2}{y_4^3}\right)}.
\end{split}
\label{x->y}
\end{equation}
where $X_1$ and $X_2$ are as in L\textsc{emma} \ref{XY}.
Furthermore, 
the morphism given by (\ref{x->y})
\begin{eqnarray*}
\ovl\sigma:&\cV_1\longrightarrow &\cV_2,\\
&(x_0,\cd,x_5)\mapsto &(y_0,\cd,y_5).
\end{eqnarray*}
is a bi-positive birational isomorphism, that is, 
there exists the inverse birational 
isomorphism $\ovl\sigma^{-1}$ and it is also positive:
\begin{eqnarray*}
&&x_0=\frac{Y_{\ovl 1}}{y_0},\q
x_1=\frac{Y_{\ovl 2}}{y_0},\q
x_2=\frac{Y_{\ovl 3}}{y_0},\q
x_4=\frac{y_2y_4Y_{\ovl 1}}{y_0Y_{\ovl 3}},\\
&&x_3=\frac{PY_{\ovl 1}}{y_0^2Y_{\ovl 2}},\q
x_5=\frac{{y_5}Y_{\ovl 1}\,\left( 1 + \frac{{y_1}}{{y_0}} 
+ \frac{{y_3}\,{y_5}}{{y_0}} 
+ \frac{{y_1}\,{y_3}\,{y_5}}{{{y_0}}^2} + 
\frac{{y_1}\,{{y_3}}^2\,{y_5}}{{y_0}\,{{y_4}}^3} \right) }
{y_0^2{x_1}\,{x_3}}
\end{eqnarray*}
where $Y_{\ovl 1}, Y_{\ovl 2}, Y_{\ovl 3}$ are 
as in L\textsc{emma} \ref{XY} and 
\begin{eqnarray*}
&&\hspace{-10pt}
P={y_0} + {y_1} + \frac{{y_0}\,{y_1}\,{y_5}}{{{y_2}}^3} + 
  2\,{y_3}\,{y_5} + \frac{2\,{y_0}\,{y_3}\,{y_5}}{{y_1}} + 
  \frac{{y_0}\,{{y_3}}^2\,{y_5}}{{{y_4}}^3} + 
  \frac{2\,{y_1}\,{{y_3}}^2\,{y_5}}{{{y_4}}^3}\\
&& + \frac{3\,{y_0}\,{y_3}\,{y_5}}{{y_2}\,{y_4}}
+\frac{3\,{y_1}\,{y_3}\,{y_5}}{{y_2}\,{y_4}}
+ \frac{3\,{y_0}\,{y_4}\,{y_5}}{{{y_2}}^2} + 
  \frac{3\,{y_0}\,{{y_4}}^2\,{y_5}}{{y_1}\,{y_2}} +
 \frac{{y_0}\,{{y_4}}^3\,{y_5}}{{{y_1}}^2} + 
  \frac{{y_1}\,{y_3}\,{{y_5}}^2}{{{y_2}}^3} \\
&&
+ \frac{3\,{{y_3}}^2\,{{y_5}}^2}{{y_1}} + 
  \frac{{y_1}\,{{y_3}}^4\,{{y_5}}^2}{{{y_4}}^6} 
+\frac{3\,{y_1}\,{{y_3}}^3\,{{y_5}}^2}{{y_2}\,{{y_4}}^4}
+ \frac{3\,{{y_3}}^3\,{{y_5}}^2}{{{y_4}}^3} + 
  \frac{3\,{y_1}\,{{y_3}}^2\,{{y_5}}^2}{{{y_2}}^2\,{{y_4}}^2} + 
  \frac{6\,{{y_3}}^2\,{{y_5}}^2}{{y_2}\,{y_4}}\\
&& + 
  \frac{3\,{y_3}\,{y_4}\,{{y_5}}^2}{{{y_2}}^2} + 
  \frac{3\,{y_3}\,{{y_4}}^2\,{{y_5}}^2}{{y_1}\,{y_2}} + 
  \frac{{y_3}\,{{y_4}}^3\,{{y_5}}^2}{{{y_1}}^2}.
\end{eqnarray*}
\end{pro}
\begin{proof}
By the direct calculations, we obtain the results. 
Indeed, certain computer softwares are useful to the calculations.
\end{proof}
Here we obtain the positive birational isomorphism 
$\ovl\sigma:\cV_1\longrightarrow \cV_2$ ($v_1(x)\mapsto v_2(y)$)
 and its inverse $\ovl\sigma^{-1}$ as above.
The actions of $e_0^c$ on $v_2(y)$ 
(respectively $\gamma_0(v_2(y))$ and 
$\vep_0(v_2(y)))$ are 
induced from the ones on 
$Y_2(y_2)Y_1(y_1)Y_2(y_4)Y_1(y_3)Y_0(y_0)Y_1(y_5)$ 
as an element of the geometric crystal $\cV_2$ 
since
$e_0\bv{\ovl 2}=e_1\bv{\ovl 2}=0$. 
Now, we define the action $e_0^c$ on $v_1(x)$ by
\begin{equation}
e_0^cv_1(x)=\ovl\sigma^{-1}\circ e_0^c\circ\ovl\sigma(v_1(x))).
\label{e0}
\end{equation}
We also define $\gamma_0(v_1(x))$ and $\vep_0(v_1(x))$ by 
\begin{equation}
\gamma_0(v_1(x))=\gamma_0(\ovl\sigma(v_1(x))),\qq
\vep_0(v_1(x)):=\vep_0(\ovl\sigma(v_1(x))).
\label{wt0}
\end{equation}
\begin{thm}
Together with (\ref{e0}), (\ref{wt0}) on $\cV_1$, we obtain a
positive affine geometric crystal $\chi:=
(\cV_1,\{e_i\}_{i\in I},
\{\gamma_i\}_{i\in I},\{\vep_i\}_{i\in I})$
$(I=\{0,1,2\})$, whose explicit form is as follows:
first we have $e_i^c$, $\gamma_i$ and $\vep_i$
for $i=1,2$ from the formula (\ref{eici}), (\ref{vep-i})
and (\ref{gamma-i}).
\begin{eqnarray*}
&&\hspace{-30pt}
e_1^c(v_1(x))=v_1(x_0,\cC_1x_1,x_2,\cC_3x_3,x_4,\cC_5x_5),\,
e_2^c(v_1(x))=v_1(x_0,x_1,\cC_2x_2,x_3,\cC_4x_4,x_5),\\
&&\text{where}\\
&&\cC_1=\frac{\frac{c\,{x_0}}{{x_1}} 
+ \frac{{x_0}\,{{x_2}}^3}{{{x_1}}^2\,{x_3}} 
+\frac{{x_0}\,{{x_2}}^3\,{{x_4}}^3}{{{x_1}}^2\,
{{x_3}}^2\,{x_5}}}{\frac{{x_0}}{{x_1}} 
+ \frac{{x_0}\,{{x_2}}^3}{{{x_1}}^2\,{x_3}} + 
\frac{{x_0}\,{{x_2}}^3\,
{{x_4}}^3}{{{x_1}}^2\,{{x_3}}^2\,{x_5}}},\q
\cC_3=\frac{\frac{c\,{x_0}}{{x_1}} 
+ \frac{c\,{x_0}\,{{x_2}}^3}{{{x_1}}^2\,{x_3}} + 
\frac{{x_0}\,{{x_2}}^3\,{{x_4}}^3}{{{x_1}}^2
\,{{x_3}}^2\,{x_5}}}{
\frac{c\,{x_0}}{{x_1}} 
+ \frac{{x_0}\,{{x_2}}^3}{{{x_1}}^2\,{x_3}} + 
\frac{{x_0}\,{{x_2}}^3\,{{x_4}}^3}
{{{x_1}}^2\,{{x_3}}^2\,{x_5}}},\\
&&\cC_5=\frac{c\,\left( \frac{{x_0}}{{x_1}} + 
      \frac{{x_0}\,{{x_2}}^3}{{{x_1}}^2\,{x_3}} + 
\frac{{x_0}\,{{x_2}}^3\,{{x_4}}^3}{{{x_1}}^2
\,{{x_3}}^2\,{x_5}} \right)}{\frac{c\,{x_0}}{{x_1}} 
+ \frac{c\,{x_0}\,{{x_2}}^3}{{{x_1}}^2\,{x_3}} + 
\frac{{x_0}\,{{x_2}}^3\,{{x_4}}^3}
{{{x_1}}^2\,{{x_3}}^2\,{x_5}}},\,
\cC_2=\frac{\frac{c\,{x_1}}{{x_2}} 
+ \frac{{x_1}\,{x_3}}{{{x_2}}^2\,{x_4}}}
{\frac{{x_1}}{{x_2}} + \frac{{x_1}\,{x_3}}{{{x_2}}^2\,{x_4}}},\,
\cC_4=\frac{c\,\left( \frac{{x_1}}{{x_2}} + 
\frac{{x_1}\,{x_3}}{{{x_2}}^2\,{x_4}} \right) }{\frac{c\,{x_1}}
{{x_2}} + \frac{{x_1}\,{x_3}}{{{x_2}}^2\,{x_4}}},\\
&&\vep_1(v_1(x))={\frac{{x_0}}{{x_1}} 
+ \frac{{x_0}\,{{x_2}}^3}{{{x_1}}^2\,{x_3}} + 
\frac{{x_0}\,{{x_2}}^3\,
{{x_4}}^3}{{{x_1}}^2\,{{x_3}}^2\,{x_5}}},\q
\vep_2(v_1(x))={\frac{{x_1}}{{x_2}} 
+ \frac{{x_1}\,{x_3}}{{{x_2}}^2\,{x_4}}},\\
&&\gamma_1(v_1(x))=\frac{x_1^2x_3^2x_5^2}{x_0x_2^3x_4^3},\q
\gamma_2(v_1(x))=\frac{x_2^2x_4^2}{x_1x_3x_5}.
\end{eqnarray*}
We also have $e_0^c$, $\vep_0$ and $\gamma_0$ on $v_1(x)$
as:
\begin{eqnarray*}
&&e_0^c(v_1(x))=v_1(\frac{D}{c\cdot E}x_0,\frac{F}{c\cdot E}x_1,
\frac{G}{c\cdot E}x_2,\frac{D\cdot H}{c^2\cdot E\cdot F}x_3,
\frac{D}{c\cdot G}x_4,\frac{D}{c\cdot H}x_5),\\
&&\vep_0(v_1(x))=\frac{E}{{{x_0}}^3\,{{x_2}}^3\,{x_3}},\qq
\gamma_0(v_1(x))=\frac{x_0^2}{x_1x_3x_5},\\
&&\text{where}\\
&&D=c^2\,{{x_0}}^2\,{{x_2}}^3\,{x_3} + 
  {x_1}\,{{x_2}}^3\,{{x_3}}^2\,{x_5} + 
  c\,{x_0}\,( {x_1}\,{{x_3}}^3 + 
     3\,{x_1}\,{x_2}\,{{x_3}}^2\,{x_4}\\
&&\qq\qq\qq + 3\,{x_1}\,{{x_2}}^2\,{x_3}\,{{x_4}}^2
+ {{x_2}}^3\,\left( {{x_3}}^2 + {x_1}\,{{x_4}}^3 + 
        {x_1}\,{x_3}\,{x_5} \right) ),\\
&&E={{x_0}}^2\,{{x_2}}^3\,{x_3} 
+ {x_1}\,{{x_2}}^3\,{{x_3}}^2\,{x_5} + 
  {x_0}\,\left( {x_1}\,{{x_3}}^3 + 
     3\,{x_1}\,{x_2}\,{{x_3}}^2\,{x_4} + 
     3\,{x_1}\,{{x_2}}^2\,{x_3}\,{{x_4}}^2 \right.\\
&&\qq \left.+  {{x_2}}^3\,\left( {{x_3}}^2 + {x_1}\,{{x_4}}^3 + 
        {x_1}\,{x_3}\,{x_5} \right)  \right),\\
&&F=c\,{{x_0}}^2\,{{x_2}}^3\,{x_3} + 
  {x_1}\,{{x_2}}^3\,{{x_3}}^2\,{x_5} + 
  {x_0}\,( c\,{x_1}\,{{x_3}}^3 + 
     3\,c\,{x_1}\,{x_2}\,{{x_3}}^2\,{x_4}\\
&&\qq\qq\qq + 3\,c\,{x_1}\,{{x_2}}^2\,{x_3}\,{{x_4}}^2
+ {{x_2}}^3\,\left( {{x_3}}^2 + c\,{x_1}\,{{x_4}}^3 + 
        c\,{x_1}\,{x_3}\,{x_5} \right)),\\
&&G=c\,{{x_0}}^2\,{{x_2}}^3\,{x_3} + 
  {x_1}\,{{x_2}}^3\,{{x_3}}^2\,{x_5} + 
  {x_0}\,( {x_1}\,{{x_3}}^3 + 
     \left( 2 + c \right) \,{x_1}\,{x_2}\,{{x_3}}^2\,{x_4}\\
&&\qq\qq\qq + \left( 1 + 2\,c \right) 
\,{x_1}\,{{x_2}}^2\,{x_3}\,{{x_4}}^2
+ {{x_2}}^3\,\left( {{x_3}}^2 + c\,{x_1}\,{{x_4}}^3 + 
        c\,{x_1}\,{x_3}\,{x_5}\right)),\\
&&H=c\,{{x_0}}^2\,{{x_2}}^3\,{x_3} + 
  {x_1}\,{{x_2}}^3\,{{x_3}}^2\,{x_5} + 
  {x_0}\,( {x_1}\,{{x_3}}^3 + 
     3\,{x_1}\,{x_2}\,{{x_3}}^2\,{x_4} \\
&&\qq\qq\qq
+ 3\,{x_1}\,{{x_2}}^2\,{x_3}\,{{x_4}}^2
+ {{x_2}}^3\,\left( {{x_3}}^2 + {x_1}\,{{x_4}}^3 + 
        c\,{x_1}\,{x_3}\,{x_5} \right)).
\end{eqnarray*}
\end{thm}
\begin{proof}
Calculating directly, we know that $\chi:=
(\cV_1,\{e_i^c\}_{i\in I},
\{\gamma_i\}_{i\in I},$ $\{\vep_i\}_{i\in I})$ satisfies the 
relations in Definition \ref{def-gc} and then it is a 
geometric crystal. The positivity is immediate from 
the above formulae.
\end{proof}
Here we denote the positive structure on $\chi$ by 
$\theta:\cV_1\longrightarrow T$. Then by Corollary \ref{cor-posi}
we obtain the ultra-discretization ${\mathcal UD}(\chi,T,\theta)$, 
which is a Kashiwara's crystal. In \cite{KNO}, we show that
such crystal is isomorphic to the limit of 
certain perfect crystal for 
the Langlands dual algebra. So we present the following 
conjecture:
\begin{conj}
The crystal ${\mathcal UD}(\chi,T,\theta)$ as above 
is the limit of coherent family of 
perfect crystals of type $\TY(D,3,4)$ in \cite{KMOY}.
\end{conj}

\bibliographystyle{amsalpha}

\end{document}
